\documentclass[12pt,a4paper,reqno]{amsart}
\usepackage[centertags]{amsmath}
\usepackage{amsfonts}
\usepackage{amssymb}
\usepackage{amsthm}
\usepackage{a4}
\usepackage{enumerate}

\newtheorem{theorem}{Theorem}

\usepackage{bbm}

\newlength{\defbaselineskip}
\newcommand{\setlinespacing}[1]%
 {\setlength{\baselineskip}{#1 \defbaselineskip}}

\begin{document}

\title[On a reduction of nonlinear evolution ... ]{On a reduction of nonlinear evolution\\ and wave type equations via non-point symmetry method }

\author[I. Tsyfra]{Ivan Tsyfra}


\maketitle
\noindent AGH University of Science and Technology, Faculty of Applied Mathematics, 30 Mickiewicza Avenue, 30-059 Krakow, Poland \\[4pt]
Institute of Geophysics of NAS of Ukraine, Kyiv, Ukraine\\[4pt]
Email: tsyfra@agh.edu.pl

\bigskip

\noindent MSC classification numbers: 35A20, 35Q55, 35Q72, 82D75 \\[4pt]

\begin{abstract}
We study the symmetry reduction of nonlinear evolution and wave type
differential equations by using operators of non-point symmetry.
In our approach we use both operators of classical and conditional symmetry.
It appears that the combination of non-point and conditional symmetry enables
us to construct not only solutions but B\"acklund transformations too for
 the equation under study.We show that the method can be applied to nonevolutionary partial
differential equations.

\end{abstract}

\vspace{0,2cm}

\section{Introduction}
In recent years the symmetry method is often used for construction solutions to different mathematical physics equations \cite{BK,O}. In this paper we study the symmetry reduction of partial differential equations by using the operators of non-point symmetry. We construct the corresponding ansatz for dependent variable $u$ or its derivatives which reduces the scalar partial differential equation to a system of ordinary differential equations. After integration of reduced ordinary differential equations one obtains partial solutions of the equations under study. Application of conditional symmetry essentially widens the class of ansatz reducing partial differential equation. By using operators of Lie--B\"acklund symmetries of ordinary differential equation we can construct an ansatz for $u$ which reduces nonevolutionary partial differential equations with two independent variables to the system of ordinary differential equations and also the number of equations is smaller then the number of unknown functions. This enables us to obtain solutions depending on arbitrary functions. It is obvious that the existence of conditional symmetry does not guarantees that the solution obtained with the help of corresponding ansatz is really new that it is not an invariant solution in the classical Lie sense. We obtain the sufficient conditions that ensure that the solution found with the help of conditional symmetry operators is an invariant one in the classical sense.

\section{Non-point symmetry and reduction of nonlinear wave type\\ and evolution equations with two independent variables }
The concept of differential invariant solutions based on infinite Lie group $G$ is introduced in~\cite{Ovs}. This group is a classical symmetry group of point transformations of dependent and independent variables for the equation under study. Generally speaking, analysis similar to that in constructing differential invariant solutions enables us to obtain the ans\"atze for derivatives~$u_{x_1}$,~$u_{x_2}$ by virtue of operators of non-point symmetry \cite{Ovs,T2}.
Let us consider nonlinear differential equation
\begin{equation}\label{2}
u_{x_2x_2}=\frac{1}{{\rm e}^{u_{x_1}}-C}, \quad C={\rm const}.
\end{equation}
We search for the ansatz for the derivatives of such form
\begin{equation}\label{1}
\frac{\partial u}{\partial x_1}=R_1(x_1, x_2, u,\varphi_1(\omega) ,\varphi_2(\omega)), \quad \frac{\partial u}{\partial x_2}=R_2(x_1, x_2, u, \varphi_1(\omega), \varphi_2(\omega)) ,
\end{equation}
where $\omega=\omega(x_1, x_2, u)$. Operators of classical and conditional symmetry of the corresponding system can be used to find $R_1$, $R_2$. The corresponding system has the form
\begin{equation}\label{3}
 v^1_2=v^2_1, \quad v^2_2=\frac{1}{{\rm e}^{v^1}-C},
\end{equation}
where $v^1=u_{x_1}$,$v^2=u_{x_2}$, $v^i_k=v^i_{x_k}$, $i, k=1,2$.
 To construct ansatz of type (\ref{1}) we use the symmetry operator $\frac{1}{2}(Q-D)$ of system (\ref{3}), where
\[
D=2x_1\partial_{x_1}+x_2\partial_{x_2}+v^2\partial_{v^2}, \quad Q=\big(x_2+2Cv^2\big)\partial_{x_2}+2\partial_{v^1}-v^2\partial_{v^2}.
\]

In this case $\omega$ does not depend on $u$ because of invariance of equation (\ref{2}) with respect to one-parameter group of translation $u'=u+a$, where $a$ is a group parameter. Wit the help of the symmetry operator one can construct the ansatz
\[
v^1=\varphi_1(\omega)-\ln(x_1\varphi_2(\omega)), \quad v^2=x_1\varphi_2(\omega), \quad \omega =Cv^2+x_2
\]
reducing (\ref{3}) to the system of ordinary differential equations
\begin{equation}\label{4}
{\rm e}^{\varphi_1}\frac{{\rm d}\varphi_2}{{\rm d}\omega}=\varphi_2, \quad \varphi_2\frac{{\rm d}\varphi_1}{{\rm d}\omega}-\frac{{\rm d}\varphi_2}{{\rm d}\omega}=
\varphi_2^2.
\end{equation}
By integrating (\ref{4}) we find
\begin{equation}\label{5}
\varphi_1=\ln \frac{C_1^2\exp(2\alpha \omega)-1}{2\alpha C_1\exp(\alpha \omega)}, \quad
\varphi_2=\frac{\alpha( C_1\exp(\alpha \omega)-1)}{C_1\exp(\alpha \omega)+1},
\end{equation}
where $C_1, \alpha$ are arbitrary real constants. Then the construction of solutions of equation (\ref{2}) requires integration of overdetermined system
\[
u_{x_1}= \frac{(C_1\exp(\alpha (x_2+Cu_{x_2} ))+1)^2}{2\alpha^2t C_1\exp(\alpha (x_2+Cu_{x_2}) )},
\]
\[
u_{x_2} =\alpha t \frac{C_1\exp(\alpha (x_2+Cu_{x_2} ))-1}{C_1\exp(\alpha (x_2+Cu_{x_2} ))+1}.
\]
However, it is easy to present the solution of the equation
\begin{equation}\label{6}
w_t+ \left (\frac{w_x}{w(Cw+1)}\right )_x=0
\end{equation}
in the form
\[
\frac{Cw+1}{w}= \frac{(C_1\exp(\alpha (x+C\theta ))+1)^2}{2\alpha^2t C_1\exp(\alpha (x+C\theta ) )},
\]
\[
\theta =\alpha t \frac{C_1\exp(\alpha (x+C\theta ))-1}{C_1\exp(\alpha (x+C\theta ))+1}.
\]

Next we show that the operators of conditional symmetry of corresponding system can be used for construction the B\"acklund transformations for nonlinear wave equation
\begin{equation}\label{7}
u_{x_1x_2}=\big[1-k^2u_{x_2}^2\big]^{1/2}\sin u.
\end{equation}
Without loss of generality we can search for the operator of conditional symmetry of the corresponding system
\begin{equation}\label{8}
v^1_2+v^1_3v^2=v^2_1+v^2_3v^1,
\end{equation}
\begin{equation}\label{9}
v^2_1+v^2_3v^1=\sqrt{1-k^2(v_2)^2}\sin x_3,
\end{equation}
where $u\equiv x_3$ in the form
\begin{equation}\label{10}
Q=\partial_{x_3}+\eta^1(x_1, x_2, x_3)\partial_{v^1}+\eta^2(x_1, x_2, x_3)\partial_{v^2}.
\end{equation}
From the infinitesimal criterion of conditional invariance we obtain equalities
\begin{equation}\label{11}
{\mathop Q\limits_1}(v^1_2+v^1_3v^2-v^2_1-v^2_3v^1 )=0,
\end{equation}
\begin{equation}\label{12}
{\mathop Q\limits_1}\left (v^2_1+v^2_3v^1-\sqrt{1-k^2(v_2)^2}\sin x_3 \right )=0,
\end{equation}
where ${\mathop Q\limits_1}$ is the first prolongation of operator $Q$ \cite{O}, which should be valid on the manifold defined by (\ref{8}), (\ref{9}) and by expression $v^1_3=\eta^1$, $v^2_3=\eta^2$. Assuming that $\eta^1$ does not depend on $v^1$ we obtain from (\ref{11})
\begin{equation}\label{13}
\frac{\eta^1_3}{\sin x_3}=\frac{\eta^2_{v^2}\sqrt{1-k^2(v_2)^2}}{v^2}=C={\rm const}.
\end{equation}
It follows from (\ref{13}) that $\eta^1$ and $\eta^2$ have the form
\begin{equation}\label{14}
\eta^1=-C\cos x_3+C_1,\quad C_1={\rm const},
\end{equation}
\begin{equation}\label{15}
\eta^2=-\frac{C}{k^2}\sqrt{1-k^2(v^2)^2}+C_2, \quad C_2={\rm const}.
\end{equation}
We derive from (\ref{12}) $C_1=C_2=0$, $C^2=k^2$. Then we obtain the operator of conditional symmetry
\begin{equation}
Q=\partial_{x_3}+k\cos x_3\partial_{v^1}+k^{-1}\sqrt{1-k^2(v^2)^2}\partial_{v^2}
\end{equation}
of the system (\ref{8}), (\ref{9}) by choosing $C=-k$. Using the operator $Q$ we construct the corresponding ansatz which can be written in the following form
\begin{equation}\label{16}
u_{x_1}=\varphi_2+k\sin u, \quad u_{x_2}=k^{-1}\sin (u-\varphi_1),
\end{equation}
where $\varphi_1$, $\varphi_2$ are unknown functions depending on $x_1$, $x_2$. Substituting (\ref{16}) into (\ref{7}) we obtain the reduced system
\begin{equation}\label{17}
\varphi_{2x_2}=\sin \varphi_1, \quad \varphi_2=\varphi_{1x_1}.
\end{equation}
Thus $\varphi_1$ satisfies the sine-Gordon equation $\varphi_{1x_1x_2}=\sin \varphi_1$. Denote by $\varphi_1=w$. Then we can rewrite (\ref{16}) in the form
\begin{equation}\label{18}
u_{x_2}=k^{-1}\sin (u-w), \quad u_{x_1}=w_{x_1}+k\sin u.
\end{equation}
These B\"acklund transformations that have been obtained in \cite{BD} by another technique transform solutions of sine-Gordon equation $w_{x_1x_2}=\sin w$ to solution of equation (\ref{7}).

Further we show the application of Lie--B\"acklund symmetry for reduction of partial differential equations. Let us consider equation
\begin{equation}\label{19}
U ( x, u,{\mathop u\limits_1},{\mathop u\limits_2},\ldots ,
{\mathop u\limits_k} ) =0,
\end{equation}
where $x=(x_1,x_2,\ldots ,x_n)$, $u=u(x)\in C^k(\mathbb{R}^n,{\mathbb
R}^1)$, and
${\mathop u\limits_k}$ denotes all partial derivatives of $k$-th order and the $m$-th order ordinary differential equation of the form
\begin{equation}\label{20}
H\left (x_1, x_2, \ldots , x_n, u, \frac{\partial u}{\partial x_1}, \ldots, \frac{\partial^m u}{\partial x_1^m} \right )=0.
\end{equation}
Let
\begin{equation}\label{21}
u=F(x, C_1,\ldots,C_{m}),
\end{equation}
where $F$ is a smooth function of variables $x, C_1,\ldots ,C_{m}$,
$C_1,\ldots, C_{m}$
are arbitrary functions of parametric variables $x_2, x_3,\ldots, x_n$,
be a general solution of equation~(\ref{20}).
We also use the Lie--B\"acklund operator in the canonical form
\begin{equation}\label{22}
X=U ( x, u,{\mathop u\limits_1},{\mathop u\limits_2},\ldots ,{\mathop u\limits_k} )\partial_u.
\end{equation}
The following statement holds \cite{T8}.
\begin{theorem}\label{t}
Let equation~\eqref{20} be invariant with respect to
the Lie--B\"acklund operator $X$. Then the ansatz
\begin{equation}\label{23}
u=F(x, \varphi_1, \varphi_2,\ldots,\varphi_{m}),
\end{equation}
where $\varphi_1, \varphi_2,\ldots,\varphi_{m}$ depend on $n-1$ variables
 $x_2, x_3,\ldots, x_n$ reduces partial differential equation~\eqref{19}
to the system of~$k_1$ equations for unknown functions
$\varphi_1, \varphi_2,\ldots ,\varphi_m$ with $n-1$ independent variables
 and $k_1 \le m$.
\end{theorem}
The theorem can be easily generalized for $\varphi_1, \varphi_2,\ldots ,\varphi_{m}$ depending
on $\omega_l(x)$, $l=\overline{1, n-1}$, where $\omega_l(x)$ are some
functions of variables $x$.

Next consider nonlinear ordinary differential equation of such form
\begin{equation}\label{32}
u_{x_1x_1}+u_{x_1}^2=0.
\end{equation}
Recall, that the concepts of local theory of differential equations such as symmetry, conditional symmetry, conservation laws, Lax representations are defined by differential equalities which must be satisfied only for solutions of the equations under study. We have proved that equation~(\ref{32}) admits two Lie--B\"aclund symmetry operators
\begin{equation}\label{33}
Q_1=\frac{u_{x_1x_2}}{u_{x_1}^2}\partial_u,\quad Q_2= F(u+\ln u_{x_1})\partial_u,
\end{equation}
where $F$ is arbitrary smooth function depending of one variable. Thus the ansatz
\begin{equation}\label{34}
u=\ln \left (x_1+\varphi_1(x_2)\right )+\varphi_2(x_2)
\end{equation}
obtained by integrating of (\ref{32}) reduces the nonlinear wave equation
\begin{equation}\label{35}
u_{x_1x_2}=u_{x_1}^2F(u+\ln u_{x_1})
\end{equation}
to one ordinary differential equation
\begin{equation}\label{36}
\varphi_1'=-F(\varphi_2).
\end{equation}
This equation is integrable by quadratures for arbitrary $\varphi_2(x_2)$. Its general solution has the form
\begin{equation}\label{38}
u=\ln \left (x_1-\int F(\varphi_2(x_2){\rm d}x_2\right )+ \varphi_2(x_2),
\end{equation}
where $\varphi_2(x_2)$ is arbitrary smooth function.
So, in the framework of this approach we have constructed
solution with arbitrary function $\varphi_2(x_2)$ to nonlinear hyperbolic type partial differential equation (\ref{35})
for arbitrary functions $F$.

In addition we formulate the theorem concerning the sufficient conditions for the solution obtained by using
conditional symmetry operators to be an invariant solution in the classical
sense. Namely consider involutive family of operators
\begin{equation}\label{37}
Q_a=\xi_{aj}(x,u)\partial_{x_j}+\eta_a(x,u)\partial_u, \quad
a= 1,\dots,p .
\end{equation}
We have summation on repeated indexes. Let equation~(\ref{19}) be conditionally invariant with
respect to involutive family of operators~(\ref{37}) and
corresponding ansatz reduces
this equation to ordinary differential equation. Suppose that general solution of reduced equation depends on $C_1, C_2, \dots, C_t$ real constants.
Then the following statement holds.
\begin{theorem}\label{t2}
Let equation~\eqref{19} be invariant with respect to $s$-dimensional Lie algebra $AG_s$ and conditionally invariant with respect to involutive family of operators $\{Q_i\}$. If the system
\[
\xi_{aj} \frac{\partial u}{\partial x_j}=\eta_a(x,u)
\]
is also invariant under the algebra $AG_s$ and $s\ge t+1$, then conditionally invariant solution of equation~\eqref{19} with respect to involutive family of operators $\{Q_a\}$ is an invariant solution in the classical Lie sense.
\end{theorem}
Theorem \ref{t2} can be applied not only to conditionally invariant solutions but to solutions constructed by method of differential constraints or potential symmetry.

\section{Conclusions}

We use the operators of non-point symmetry to construct ans\"atze for dependent variable $u$ and its derivatives. It turns out that the conditional symmetry of corresponding system enables us to construct algorithmically the B\"acklund transformations to nonlinear wave equation. It should be noted
that one can construct ans\"atze for derivatives by using operators of point symmetry admitted by the initial equation but they lead to the
 invariant solutions in the classical Lie sense. Thus to obtain new results it is necessary to use operators of non-point and conditional symmetry. We obtained theorem allowing us to exclude the operators that lead to the classical invariant solutions. We showed also that the operators of Lie--B\"acklund symmetry can be used for reducing partial differential equations which are not restricted to evolution type ones.
Finally we show that the application of Theorem~\ref{t} gives the possibility of constructing the solution~(\ref{38}) defined by arbitrary function to equation~(\ref{35}). To our knowledge the inverse scattering transformation method is not applicable in this case.


\begin{thebibliography}{99}

\small

\bibitem{BK}
Bluman G.W. and Kumei S.,
Symmetries and Differential Equations,  Appl. Math. Sci., Vol.~81, Springer-Verlag, Berlin, 1989.

\bibitem{O}
Olver P.J., Applications of Lie Groups to Differential
Equations, Springer-Verlag, New York, 1986.

\bibitem{Ovs}
Ovsiannikov L.V., Group Analysis of Differential Equations, Academic Press, New York, 1982.


\bibitem{T2}
Tsyfra I., Napoli A., Messina A. and Tretynyk V. On new ways of group methods for reduction of evolution-type
equations, {\it J. Math. Anal. Appl.} {\bf 307} (2005), 724--735.

\bibitem{BD}
Dodd R.K. and Bullogh R.K. B\"acklund transformations for the sine-Gordon equations, {\it Proc. R. Soc.~A} {\bf 351} (1976), 499--523.

\bibitem{T8} Tsyfra I.M. Symmetry reduction of nonlinear differential equations, \emph{Proceedings of Institute of Mathematics, Kiev} {\bf 50} (2004), 266--270.

\end{thebibliography}
\end{document}